\documentclass[12pt]{amsart}
\usepackage{amsmath,amsfonts,amssymb,amsthm}
\usepackage{pinlabel}
\usepackage{graphicx}
\usepackage{mathtools}
\usepackage[all]{xy}
\usepackage{hyperref}
\usepackage[usenames,dvipsnames]{color}
\usepackage{color}
\xyoption{dvips}

\newtheorem{theorem}{Theorem}

\newtheorem{corollary}[theorem]{Corollary}

\theoremstyle{remark}
\newtheorem{remark}[theorem]{Remark}

\theoremstyle{example}
\theoremstyle{definition}

\newcommand{\R}{\mathbb{R}}

\numberwithin{equation}{section}

\begin{document}

\title[Lagrangian concordance is not a partial order in high dimensions]{Lagrangian concordance is not a partial order in high dimensions}

\author{Roman Golovko}

\begin{abstract}
In this short note we provide the examples of pairs of closed, connected Legendrian non-isotopic Legendrian submanifolds $(\Lambda_{-}, \Lambda_{+})$ of the $(4n+1)$-dimensional contact vector space, $n>1$, such that there exist Lagrangian concordances from $\Lambda_-$ to $\Lambda_+$ and from $\Lambda_+$ to $\Lambda_-$. This contradicts anti-symmetry of the Lagrangian concordance relation, and, in particular, implies that  Lagrangian concordances with connected Legendrian ends do not define a partial order in high dimensions.
In addition, we explain how to get the same result 
for the relation given by exact Lagrangian cobordisms with connected Legendrian ends in the $(2n+1)$-dimensional contact vector space, $n>1$.
\end{abstract}

\address{Faculty of Mathematics and Physics, Charles University, Sokolovsk\'{a} 83, 18000 Praha 8, Czech Republic}
\email{golovko@karlin.mff.cuni.cz}
\date{\today}
\thanks{}


\maketitle

\section{Introduction and main result}

Lagrangian cobordisms between Legendrian submanifolds naturally arise in Symplectic Field Theory of Eliashberg, Givental and Hofer \cite{SFT00}. They have been actively studied for more than two decades, for the basics of this theory we refer to \cite{LagCobGenFam, ChantraineLagrConc,  FlhomLagConc, Lagracobbasics}.  
In this paper we concentrate on the most basic Lagrangian cobordisms, namely Lagrangian concordances with connected Legendrian ends, in the symplectization of the contact vector space $\R^{2n+1}_{st}:=(\R^{2n+1}, dz-ydx)$. 
Lagrangian concordances with connected Legendrian ends define a relation on closed, connected Legendrian submanifolds, i.e. we say that $\Lambda_-$ is related to $\Lambda_+$ and write $\Lambda_- <_{L} \Lambda_+$ if there is a Lagrangian concordance $L$ from $\Lambda_-$ to $\Lambda_+$. This relation was studied in \cite{ChantraineLagrConc, ObsLagConc}. It  is reflexive and transitive. 
By the result of Chantraine it follows that this relation is not  symmetric \cite{LagrConc_not_sym15}, see also \cite{ObsLagConc}. 
For a while there was an open question of whether this relation defines a partial order, see \cite{ ChantraineLagrConc, LagUppBound, LagZigCob}.
Recently, Dimitroglou Rizell and the author in \cite{DimitroglouRizellGolovko24} constructed 
a series of examples of Lagrangian concordances between Legendrian knots that contradict the anti-symmetry property. In other words, in \cite{DimitroglouRizellGolovko24} Dimitroglou Rizell and the author constructed the series of pairs $(\Lambda_-, \Lambda_+)$ of Legendrian non-isotopic Legendrian knots such that there are Lagrangian concordances from $\Lambda_-$ to $\Lambda_+$ and from $\Lambda_+$ to $\Lambda_-$. 

In this note, we construct the series of examples of Lagrangian concordances with connected Legendrian  ends that contradict anti-symmetry in high dimensions.
More precisely, we prove the following:
\begin{theorem}
\label{th:main_ans}
There is a pair of closed, connected Legendrian non-isotopic Legendrian submanifolds $(\Lambda_-, \Lambda_+)$ of $\R^{4n+1}_{st}$, $n>1$, such that there exist two Lagrangian concordances $L_{\pm}$ from $\Lambda_-$ to $\Lambda_+$ and $L_{\mp}$ from $\Lambda_+$ to $\Lambda_-$ in the symplectization of $\R^{4n+1}_{st}$, $n>1$.
\end{theorem} 
Theorem~\ref{th:main_ans} leads to the following corollary:
\begin{corollary}
Lagrangian concordances with connected Legendrian ends do not define a partial order on closed, connected Legendrian submanifolds of $\R_{st}^{4n+1}$ for $n>1$.
\end{corollary}
In addition, in Section \ref{Lagcobordisms} we explain how to get the same result for the relation given by exact Lagrangian cobordisms with connected Legendrian ends in $\R^{2n+1}_{st}$, $n>1$.

\section{Proof of Theorem~\ref{th:main_ans}}
First we recall that from the result of Murphy  \cite[Theorem 1.2]{LooseLegendrians} and the discussion right after  it follows  that every formal Legendrian isotopy class in $\R^{2k+1}_{st}$, $k>1$, contains loose Legendrian embeddings, and every two formally isotopic loose  Legendrian
embeddings are Legendrian isotopic.

Then we assume that $k=2n$, $n>1$. Following the result of Ekholm-Etnyre-Sullivan \cite[Proposition 3.2]{EkholmEtnyreSullivan05} observe that for a Legendrian submanifold $\Lambda$ of $\R^{4n+1}_{st}$, the Thurston-Bennequin number $tb(\Lambda)$ is a topological invariant given by 
$$tb(\Lambda)=(-1)^{n+1}\frac{\chi(\Lambda)}{2}.$$

Now we take a closed, stably parallelizable, simply connected $\Lambda$ and 
two Legendrian embeddings $i_-:\Lambda\to \R^{4n+1}_{st}$ and $i_+:\Lambda\to \R^{4n+1}_{st}$ such that $r(i_-(\Lambda))=r(i_+(\Lambda))$ and $\Lambda_-:=i_-(\Lambda)$, $\Lambda_+:=i_+(\Lambda)$ are loose Legendrian non-isotopic Legendrian submanifolds. 
From \cite[Proposition A.4 (c)]{LooseLegendrians} it follows that for a fixed rotation class, there is exactly one such couple $\Lambda_-$, $\Lambda_+$ up to Legendrian isotopy.
\begin{remark}
Note that, in particular, for our construction we can take $\Lambda=S^{2n}$.
\end{remark}
\begin{remark}
Following \cite[Proposition 3.2]{EkholmEtnyreSullivan05} we see that 
$tb(\Lambda_-)=tb(\Lambda_+)$. By the original assumption, the  rotation classes of  $\Lambda_-$ and of  $\Lambda_+$ coincide. Since $\Lambda_{-}$, $\Lambda_{+}$ are loose, they have acyclic Legendrian contact homology DGAs. Hence, we conclude that $\Lambda_-$ and $\Lambda_+$ have the same Legendrian invariants.
\end{remark}

From the classical result of Wu \cite{Wu58} it follows that every two embeddings of a connected $2n$-dimensional manifold into $\R^{4n+1}$, $n\geq 1$, are smoothly isotopic. Hence, $\Lambda_-$ is smoothly isotopic to $\Lambda_+$. Thus, in particular,  there is a smooth concordance  from $\Lambda_-$ to $\Lambda_+$ and we denote the corresponding  smooth concordance by $L^{C^{\infty}}_{\pm}$, and there is a smooth concordance from $\Lambda_+$  to $\Lambda_-$ and we denote it  by $L^{C^{\infty}}_{\mp}$.

For simplicity, from now on let us assume that $\Lambda=S^{2n}$. Clearly $\Lambda_-, \Lambda_+ \simeq S^{2n}$ are closed, connected, loose, and the complexified tangent bundle of a concordance over $S^{2n}$ is trivial. 
Now recall that in high dimensions Eliashberg and Murphy  establish an h-principle for exact Lagrangian embeddings with loose concave Legendrian
ends \cite[Theorem 2.2]{EliashbergMurphy13}.

\begin{remark}
\label{El_Mur_cobordisms}
Even though \cite[Theorem 2.2]{EliashbergMurphy13} is written for exact Lagrangian caps with loose concave Legendrian ends, the proof of \cite[Theorem 2.2]{EliashbergMurphy13} can be adapted to hold for exact Lagrangian cobordisms with loose concave Legendrian ends and possibly non-trivial convex Legendrian ends. 
Observe that this adaptation is possible, since all the homotopies and isotopies  from \cite[Theorem 2.2]{EliashbergMurphy13} and \cite[Theorem 2.3]{EliashbergMurphy13}, that is the statement from which \cite[Theorem 2.2]{EliashbergMurphy13} is deduced, are compactly supported.

\end{remark}

By applying the $h$-principle of Eliashberg--Murphy that we discussed in Remark \ref{El_Mur_cobordisms}
to concordances $L^{C^{\infty}}_{\pm}$, $L^{C^{\infty}}_{\mp}$ with loose ends over $\Lambda_-$, $\Lambda_+$ we get Lagrangian concordances $L_{\pm}$ from $\Lambda_-$ to $\Lambda_+$ and $L_{\mp}$ from $\Lambda_+$ to $\Lambda_-$.  Since $\Lambda_-$ is not Legendrian isotopic to $\Lambda_+$, this leads to the fact that Lagrangian concordances with connected Legendrian ends do not define an anti-symmetric relation. In particular, it implies that Lagrangian concordances with connected Legendrian ends do not define a partial order for closed, connected Legendrian submanifolds of $\R^{4n+1}_{st}$, where $n>1$. This finishes the proof. 
\begin{remark}
Note that the the proof of Theorem \ref{th:main_ans} heavily relies on the description  of Legendrian isotopy classes of loose Legendrian submanifolds in $\R^{2k+1}_{st}$ written by Murphy in  \cite[Proposition A.4]{LooseLegendrians}.  
For $k>2$ even, a given stably parallelizable, simply connected $k$-manifold $\Lambda$, and a given rotation class $r$, this classification provides  two Legendrian non-isotopic Legendrian embeddings of $\Lambda$ realizing $r$.
The fact that there are two Legendrian non-isotopic loose Legendrian embeddings of $\Lambda$ realizing the same pair of classical invariants $(r,tb)$ is crucial for the proof. Now note that for odd $k>1$,  Murphy in  \cite[Proposition A.4]{LooseLegendrians}
proved that if two formal Legendrian embeddings have the same
Thurston-Bennequin number and rotation class, then they are formally Legendrian isotopic. This, together with \cite[Theorem 1.2]{LooseLegendrians} and the discussion right after  it, implies that for a given stably parallelizable $k$-manifold $\Lambda$ and a pair 
of classical Legendrian invariants $(r, tb)$, there is exactly one loose Legendrian  embedding realizing this pair. This, in particular, implies that the method we used for the proof of Theorem \ref{th:main_ans} for even $k$ can not be easily adjusted to the case when $k$ is odd.  
\end{remark}

\section{Exact Lagrangian cobordisms}
\label{Lagcobordisms}
If instead of the class of Lagrangian concordances with connected Legendrian ends we considers the class of  exact Lagrangian cobordisms with connected Legendrian ends in $\R^{2n+1}_{st}$, $n>1$, and define the corresponding relation, then one can construct examples of Lagrangian cobordisms with non-trivial topology that contradict  the anti-symmetry. 

\begin{figure}[h!]
\begin{center}
\includegraphics[width=400pt]{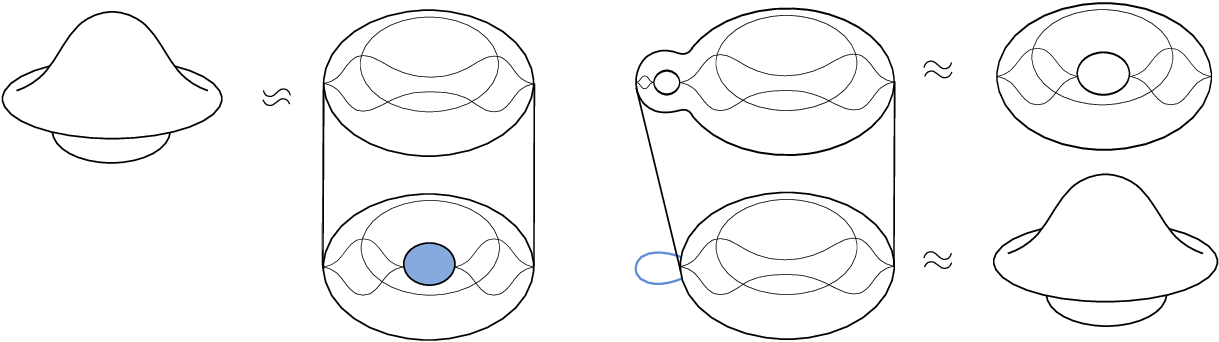}
\caption{The pair of exact Lagrangian cobordisms $L_{T^2}^{S^2}$ (left) and $L_{S^2}^{T^2}$ (right) from \cite[Section 2.3]{DimitroglouRizellGolovko14}.}\label{nonant}
\end{center}
\end{figure}

For that, we consider the example from \cite[Section 2.3]{DimitroglouRizellGolovko14}, where Dimitroglou Rizell and the author constructed the non-cylindrical exact Lagrangian endocobordism which consists of the following concatenated pieces:
\begin{itemize}
\item exact Lagrangian cobordism $L_{S^2}^{T^2}$  from the standard Legendrian $2$-sphere to the Legendrian torus which is a front spun of $tb=-1$ Legendrian unknot, and
\item exact Lagrangian cobordism $L_{T^2}^{S^2}$ from  the Legendrian torus which is a front spun of $tb=-1$ Legendrian unknot to the standard Legendrian $2$-sphere. 
\end{itemize}
These cobordisms appear in Figure \ref{nonant}.

This pair of cobordisms $L_{S^2}^{T^2}$, $L_{T^2}^{S^2}$ contradicts the anti-symmetry property by the obvious topological reason, i.e. by the fact that $S^2$ is not diffeomorphic to $T^2$. This example can easily be extended to high dimensions by using the front spinning construction (or the spherical spinning construction) applied to  $L^{S^2}_{T^2}$ and $L_{S^2}^{T^2}$ as explained in \cite{LagrCobG13, FrontSpunG14}. This 
way we obtain the exact Lagrangian cobordisms from the Legendrian $\Lambda$ diffeomorphic to $S^{i_1}\times \dots \times S^{i_k}\times S^2$ to the Legendrian $\Lambda'$ diffeomorphic to  $S^{i_1}\times \dots \times S^{i_k}\times T^2$ and from $\Lambda'$ to $\Lambda$, which implies that the relation defined by exact Lagrangian cobordisms with connected Legendrian ends for closed, connected Legendrian submanifolds of $\R^{2n+1}_{st}$ is not anti-symmetric for all $n>1$, and hence it is not a partial order for all $n>1$.

\section*{Acknowledgements}
The author would like to thank Baptiste Chantraine, Georgios Dimitroglou Rizell and Oleg Lazarev for the very helpful discussions. 
The author is supported by the 4EU+ academic mini-grant MA/4EU+/2024/F3/01.


\begin{thebibliography}{}

\bibitem{LagCobGenFam}
F. Bourgeois, J. Sabloff and L. Traynor, {\em Lagrangian cobordisms via generating families:
construction and geography}, Algebraic and Geometric Topology 15 (2015) 2439--2477.

\bibitem{ChantraineLagrConc}
B. Chantraine, {\em On Lagrangian concordance of Legendrian knots}, Algebr. Geom. Topol. 10 (2010), 63--85.

\bibitem{LagrConc_not_sym15}
B. Chantraine, {\em Lagrangian concordance is not a symmetric relation}, Quantum Topol. 6
(2015), no. 3, 451--474.

\bibitem{FlThLagCob}
B. Chantraine, G. Dimitroglou Rizell, P. Ghiggini and R. Golovko, {\em
Floer theory for Lagrangian cobordisms}, Journal of Differential Geometry, Volume 114, Number 3 (2020), 393--465.

\bibitem{FlhomLagConc}
B. Chantraine, G. Dimitroglou Rizell, P. Ghiggini and R. Golovko, {\em Floer homology and Lagrangian concordance}, 
Proceedings of the Gökova Geometry-Topology Conference 2014, (2015), 76--113.

\bibitem{ObsLagConc}
C. Cornwell, L. Ng, and S. Sivek, {\em Obstructions to Lagrangian concordance}, Algebr.
Geom. Topol. 16 (2016), no. 2, 797--824.

\bibitem{DimitroglouRizellGolovko14}
G. Dimitroglou Rizell and R. Golovko, {\em On homological rigidity and flexibility of exact Lagrangian endocobordisms}, 
International Journal of  Mathematics, Vol. 25, No. 10, 1450098 (2014), 24 pages.

\bibitem{DimitroglouRizellGolovko24}
G. Dimitroglou Rizell and R. Golovko, {\em Instability of Legendrian knottedness, and non-regular Lagrangian concordances of knots}, preprint 2024, available at arXiv:2409.00290. 


\bibitem{EkholmEtnyreSullivan05}
T. Ekholm, J. Etnyre, and M. Sullivan, {\em Non-isotopic legendrian
submanifolds in $\R^{2n+1}$}, Journal of Differential Geometry, volume 71(1):85 --128, 2005.

\bibitem{Lagracobbasics}
T. Ekholm, K. Honda and T. K\'{a}lm\'{a}n, {\em Legendrian knots and exact Lagrangian cobordisms},
J. Eur. Math. Soc. 18 (2016), no. 11, pp. 2627--2689.

\bibitem{EliashbergMurphy13}
Y. Eliashberg and  E. Murphy,
{\em Lagrangian caps}, Geometric and Functional Analysis, Volume 23, 1483--1514, 2013. 

\bibitem{SFT00}
Y Eliashberg, A Givental and H Hofer, {\em Introduction to
symplectic field theory}, Geom. Funct. Anal. (2000), 560--673.

\bibitem{LagrCobG13}
R. Golovko, {\em A note on Lagrangian cobordisms between Legendrian submanifolds of $\R^{2n+1}$}, Pacific Journal of Mathematics, Vol. 261 (2013), No. 1, 101--116.

\bibitem{FrontSpunG14}
R. Golovko, {\em A note on the front spinning construction}, Bulletin of the London Mathematical Society, 46 (2014), no. 2, 258--268.

\bibitem{LooseLegendrians}
E. Murphy, {\em Loose Legendrian embeddings in high dimensional
contact manifolds}, preprint 2012, available at
\texttt{arXiv:1201.2245}.

\bibitem{LagUppBound}
J. Sabloff, D. Shea Vela-Vick, C.-M. Wong,
{\em Upper bounds for the Lagrangian cobordism relation on Legendrian links}, preprint 2021. 
Available at arXiv:22105.02390.

\bibitem{LagZigCob}
J. Sabloff, D. Shea Vela-Vick, C.-M. Wong, A. Wu
{\em Lagrangian zigzag cobordisms}, preprint 2023. Available at arXiv:2308.02057.

\bibitem{Wu58}
W. T. Wu, {\em On the realization of complexes in a Euclidean space}, I, Sci Sinica 7 (1958), 251--297;
II, Sci Sinica 7 (1958), 365--387; III, Sci Sinica 8 (1959), 133--150.

\end{thebibliography}
\end{document}